\documentclass[11pt,a4paper,reqno]{amsart}
\usepackage{epsfig,amssymb,latexsym,amscd,amsmath,amsthm}
\usepackage[utf8]{inputenc}

\usepackage[all,2cell]{xy}
\xyoption{v2}

\theoremstyle{plain}

\newtheorem{teo}{Theorem}[section]
\newtheorem{theo}[teo]{Theorem}

\newtheorem{lema}[teo]{Lemma}

\theoremstyle{remark}

\theoremstyle{definition}

\newtheorem{defi}[teo]{Definition}

\renewcommand{\k}{{\mathbb K}}  
\newcommand{\K}{{\mathbb K}}

\newcommand{\R}{{\mathbb R}}    
\newcommand{\C}{{\mathbb C}}

\begin{document}

\title[A mathematician of the XX$\text{th}$ Century]{Gerhard Hochschild (1915/2010) \\ A mathematician of the XX$\text{th}$ Century}

\author{Walter Ferrer Santos}
\address{Facultad de Ciencias\\Universidad de la Rep\'ublica\\ Igu\'a 4225\\11400 Montevideo\\Uruguay.}
\thanks{The  author would like to thank Anii, Csic-UdelaR,  Conycit-MEC, Uruguay and Mathamsud Project for financial assistance. Moreover the author wants to warmly thank Hochschild's family for letting him share some parts of Gerhard's files and for providing some translations. Moreover, we have also 
used regularly along this note, the material presented in {\em Gerhard Hochschild (1915--2010)} \cite{kn:memorial} and we want to thank all the authors of the collaborations therein.}
\email{wrferrer@cmat.edu.uy}

\begin{abstract} Gerhard Hochschild's contribution to the development of mathematics in the XX century is succinctly surveyed. We start with a personal and mathematical biography, and then consider with certain detail his contributions to algebraic groups and Hopf algebras.  
\end{abstract}

\maketitle

\section{The life, times and mathematics of Gerhard Hochschild}
\label{section:life}

\subsection{Berlin and South Africa}
Gerhard Paul Hochschild was born on April 29th, 1915 in Berlin of a middle class Jewish family and died on July 8th,  2010 in El Cerrito where he lived after he moved to take a position as professor at the University of California at Berkeley in 1958. 

His father Heiner was an engineer working as a patent attonery and in search of safety sent Gerhard and his older brother Ulrich, to Cape Town in May of 1933. The boys were some of the many germans escaping from the Nazis that were taking over their native country.

The first 18 years of his life in Germany were not uneventful. In 1924  his mother Lilli was diagnosed with a lung ailment and sent --together with his younger son Gerhard that was then nine years old-- to a sanatorium in the Alps, near Davos. Later in life he commented to the author that reading in {\em Der Zauberberg} (The Magic Mountain) by Thomas Mann, about 
Hans Castorp --Mann's main character in the book-- he evoked his own personal experiences. Castorp was transported away from his ordered and organized family life to pay a visit to his cousin interned also in a sanatorium in Davos. Once he was there, Castorp was diagnosed with tuberculosis, spent seven years interned, and was only able to leave the place after volunteering for the army at the beggining of the first world war. 

While in the sanatorium his mother started to descend into mental illness and the boy was her only support. In 1926 she was transferred to a mental asylum in Germany, where she was later murdered by the Nazis as part of the infamous ``Final Solution'' program\footnote{For the Nazis, Jewish mentally-ill patients were unique among victims in that they embodied both ``hazardous genes'' and ``racial toxins'', see \cite{kn:mentallyill}.}.
\medskip 
 

Gerhard returned to Berlin to
continue his education entering the gymnasium where, even though he did not enjoy formal education --a feeling that he continued to hold all his life, he liked his courses in physics and mathematics. In a letter sent later to his mathematics teacher Dr. Flatow that was sent in 1937 after he finished a MSc. in mathematics at the University of Cape Town, he wrote: ``{\em  You probably won't remember me since it has been five years since you had me as a student.
I want to let you know how my life has gone.
I have made a choice for which I want to thank you.
You kindly advised me ... so I owe you a report.
I decided in the end for pure mathematics and not physics, although in my first two years I was involved with physics and applied mathematics...The reason for this was because I was more and more interested in mathematics, and, I came to see physics as only one field of application...I still remember with pleasure our hours of mathematics in school, and am so grateful to you for interesting me in mathematics.
Greetings from an old student''} \footnote{See Schwartz's contribtion in \cite{kn:memorial}.}.

At that time as a young boy in Germany, he started to cultivate  passions for photography and hiking in the mountains,  passions that lasted all his life.

Even though, later in life he never made a fuzz about the subject,  at the end of his stay in Berlin he suffered some of the ignominies imposed to jews by the Nazi regime.  He mentioned a particular incident to the author, in which he went to meet his young friend Eva, was accosted by Nazi thugs and one of them said while exhibiting a gun: ``{\em A german girl like you should not be with a dirty jew''}. They managed to walk away unharmed.

In South Africa, and because Hitler made it impossible to send money out of the country, the Hochschild boys were forced to earn their own way, and 
Gerhard found employment in a photographer's shop as his assistant. By then, Gerhard had started to frequent a circle of leftist intellectuals and artists that meet on a regular basis in Cape Town; later in life he frequently mentioned how well and at ease he felt within this group of free thinkers and intellectuals. It was probably at that time that he acquired the leftist and anarchist ideals that accompanied him all along his life.\footnote{When Hochschild found out about the author's involvement with some human right causes in South America while  a graduate student at Berkeley under his direction, he offered a monetary donation to :``{\em whatever organization you find suitable''.}}

His registration at the University of Cape Town in the B.Sc. program starting in January 1934, was made possible by a small grant received from a foundation established by one of his relatives: Berthold Hochshild, a german industrialist that emigrated to the USA at the end of XIX century to establish a company called American Metal Co. that flourished in the economic boom that followed the Civil War. 

He finished a B.Sc. in Physics and Mathematics in 1936 and a M.Sc.  in Mathematics in 1937. His student records list the 
following courses: {\bf 1934};  Applied mathematics I, Chemistry I, Physics I, Pure mathematics I; {\bf 1935}; Applied mathematics II, Economics I, Physics II, Pure mathematics II; {\bf 1936}; Applied mathematics III, French I, Pure mathematics III; {\bf 1937}; Applied mathematics: Relativity, Hydrodynamics, Tensor methods in dynamics; Pure mathematics: Complex analysis, Elliptic functions, Harmonic analysis, Differential geometry, Differential equations, Elementary algebra, Elementary theory of functions of a real variable.    

After finishing his M. Sc. he worked as a Junior Lecturer at the University of Cape Town during the 1937/38 academic year.

Stanley Skewes\footnote{Stanley Skewes (1899/1988) a student of Litlewood at Cambridge University is known for his discovery of the Skewes number in 1933. This is an upper bound for the smallest number $x$ for which $\pi(x) > li(x)$ where $\pi$ is the prime counting function and $li$
is the logarithmic integral; Skewes proved that the number thus defined is smaller than  $10^{10^{10^{963}}}$.}, then a lecturer at Cape Town University was Gerard's advisor and supporter. In the letter he sent to Princeton recommending Hochschild he wrote:  
``{\em I know of no reason why he should be unsuitable to enter Princeton University--Graduate College. He is a good student and a very promising mathematician --he was first choice for the post graduate scholarship at the end of 1937. I myself recommended him to go to your institution --the average University can do little more for him, except possible to set him up in research, an occupation for which I believe him to be eminently suited''.} 

In his application form to Princeton, besides attaching recommentation letters from Brown, Crawford and Skewes, Gerhard mentions a non published paper --now lost-- that he had written on: ``{\em the application of the $\delta$ tensors to the theory of determinants''}. 

He was admitted in the PhD program at Princeton University 
and in the summer 1938, he sailed from Cape Town to 
New York where his father was then living.

\subsection{Princeton and Aberdeen proving grounds}

Hochschild registered at Princeton, starting September 8, 1938\footnote{Here and in other parts of the manuscript the author used the material provided to Hochschild's family by the staff of Seeley G. Mudd Manuscript Library, Princeton University.} and  
for one of the leading algebraists of his generation his 
choice of courses in the PhD program may seem rather peculiar to the reader accustomed to the fashions of contemporary overspecialized mathematical education.

\begin{description}
\item[1938/39] Calculus of variations, I.A.S. (Mayer); Elementary theory of functions of a real variable (Bohnenblust); Continuous groups (Eisenhart); Advanced theory of functions of a real variable (Bochner); The theory of relativity (Robertson); Continuous groups (Eisenhart).
\item[1939/40] Aplications of the theory of functions of a complex variable (Strodt);
Riemannian geometry (Eisenhart); Topological groups (Chevalley); Algebraic geometry (Chevalley); Applications of analysis to geometry (Bochner);
Riemannian geometry (Eisenhart).
\item[1940/41] Applications of analysis to geometry (Bochner); Probability and ergodic theory, I.A.S. (Halmos and Ambrose); Differential equations (Chevalley); Research and work on dissertation under the direction of Chevalley
–two semesters–; Advanced theory of functions of a real variable (Bochner), Ergodic theory, I.A.S. (von Neumann).
\end{description}
 
During the years 1938/40 he enjoyed a Porter Scholarship from the University of Cape Town, and in his last year 1940/41 he was awarded in the first semester a Research Assistantship in mathematics and in the second semester a position of Part Time Instructor and Assistant. 

Hochschild was the first of Chevalley's\footnote{Claude Chevalley was a native of South Africa (1909--1984) --his father was a french diplomat-- and was at Princeton Inst. in the year 38/39, at Princeton Univ. from 1940/48 and at Columbia from 1949/55. He was one of the founding members of Bourbaki.} students at 
Princeton University, and started to work under his direction around 1939, when he gave him to study some of the first Bourbaki manuscripts.


An interesting story of his days at Princeton is the following: during the first lectures of Chevalley's course on 
Differential Equations in 1940, the room was packed with people curious to know what he had to say on this subject but at the end of the course only three persons remained: Hochschild, von Neumann and Weyl.

In the internet project: ``The Princeton Mathematics community in the 1930s'' Robert Hooke describes Chevalley as playing an ``{\em endless game of Go''} while at Princeton. Gerhard developed then, a lasting passion for the game acquiring approximately the level of a 7th. kyu. He used to play it whenever he had the oportunity, and besides Chevalley some of his rivals were: P. Erd\"os, D. Goldschmidt, N. Steenrod, etc.

While still a graduate student at Princeton, Hochschild submitted his first paper for publication. It contained the main results of his thesis, was entitled {\em Semi-simple algebras and generalized derivations} and appeared in print in 1942 shortly after he was drafted into the army. In the Introduction to this paper, see \cite{kn:hochss}, Hochschild says that it deals with: {\em ``the study of the behaviour of Lie
algebras and associative algebras with respect to derivations''} and continues: {\em ``These `generalized derivations' ... were found to be
significant for the structure of an algebra. In fact we shall obtain a
characterization of semisimple Lie algebras and semisimple associative algebras,
in terms of these generalized derivations.''}

Gerhard's dissertation committee, cochaired by Chevalley and Lefschetz, that met to hear the thesis in April 24, 1941 reports :
``{\em The thesis deals with certain important problems in Lie algebras and related questions in associative algebras. It contains in particular a highly interesting characterization of semisimple Lie algebras in terms of the operation of formal derivation. The thesis is highly worthy of publication as it contains many new results in addition to those indicated above. Furthermore Hochschild set the problem himself and also 
did the research in an essentially independent way''}.

After defending his thesis, he was appointed as a Part time Instructor and Research assistant at Princeton University, for the academic year 1941/42 starting in September, but in November he was drafted into the US Army. He was first 
stationed at Ft. Sill in Oklahoma and along with other non–citizen soldiers in his unit, he was taken to the local county courthouse for naturalization in June 1942.\footnote{See the contribution by Magid in \cite{kn:memorial}.}  
He spent the rest of his time at the Aberdeen proving grounds where he was put to work in the Mathematics Section whose director was Oswald Veblen the famous topologist from Princeton. Veblen gathered in that unit, an important number of recruits, some of them become later well known mathematicians:  Federer, Kelley, Morrey, Morse, etc.

Regarding that period of his life, in a letter sent to the 
Notices of the AMS in July 2009 and that concerns calculus teaching, Hochschild says: ``{\em Memories dating back to 1942 are of spending hundreds of days calculating military firing tables with the help of a Friden desk calculator''.}

The second and third of his papers, dealing with what later was to be called Hochschild cohomology, \cite{kn:hochcohom1} and \cite{kn:hochcohom2}
list his address as ``Aberdeen proving grounds''. 

About that period of his life, his long time collaborator G. D. Mostow says, see \cite{kn:memorial}: ``{\em One cannot write about Gerhard comprehensively without mentioning his charisma. Some of his charisma resulted from his colorful criticism of the hypocrisies that abound in all large organizations.  In the army, even though he was a recent immigrant to the United States, he impressed his fellow soldiers with the virtuosity of his profanity. I learned this from the famous geometric measure theory mathematician Herbert Federer, who served in Gerhard's unit at Aberdeen Proving Grounds''.}

After the war in mid 1945, he left the Army to take a part time position for a semester --November 1945 to June 1946-- as an Instructor at Princeton University. 

\subsection{Harvard and Urbana}

Gerhard was a Benjamin Pierce Instructor at Harvard University during the academic years 1946/48 and in September of that year, took a position at the University of Illinois at Urbana--Champaign. At Urbana he raised quickly and became a full professor in 1952. 

While tenured in Urbana, he spent the academic year 1951/52 visiting Yale University, 1955/56 visiting the University of  California at Berkeley and 1956/57 as a member of the Institute of Advanced Study at Princeton with a Guggenheim fellowship. 

Dan Mostow was also a member of the Institute at that time and their very long and fruitful collaboration started then; they wrote 17 papers together. Other Institute members that 
year were: M. Auslander, R. Bott, C. Curtis, J. Leray, I. Kaplansky, A. Rosenberg, J.-P. Serre, T. Tannaka, all of whom devolved long term personal and mathematical relationships with Gerhard. 
 

At the beggining of the 1950s, he started to impel the application to other parts of mathematics of the cohomological methods already developed and succesfully aplied to the theory of associative and Lie algebras. 

Particularly remarkable are his papers \cite{kn:lcft}, \cite{kn:arl}, \cite{kn:ccft} where he applied cohomological methods to Class field theory. In \cite{kn:lcft} he was the first to apply these methods, to the local theory. The importance  of cohomology for the global theory  was clarified in the paper \cite{kn:ccft} that was written with Nakayama. This work was explained and further developed in the famous Artin--Tate seminar, whose notes were for a long time the basic reference for the modern presentation of the theory\footnote{For this part and the next comments, see the contribution by J. Tate --with a letter from Serre-- in \cite{kn:memorial}.}.

J.-P. Serre published a short note \cite{kn:cdg}, where he sketched a proof of the construction of a spectral sequence relating the cohomolgies of $G$, $K$ and $G/K$, whenever $K$ is a normal subgroup of the finite group $G$. He used the Cartan--Leray spectral sequence associated to the action of $G/K$ on $E_G/K$ where $E_G$ is a universal bundle for $G$. After reading the note and producing a direct proof using explicit calculations with resolutions, Gerhard wrote Serre in January 1951 proposing a joint publication which included both methods. 

The first paper product of that collaboration, see \cite{kn:cge}, was written jointly while Hochschild spend one year at Yale University substituting Jacobson during his sabbatical and Serre was visiting Princeton Institute in January/February 1952 (see J.--P. Serre's letter in the contribution by J. Tate in \cite{kn:memorial}: {\em Memories of Hochschild, with a letter from Serre}). In this paper, the convergence of the so called Hochschild--Serre spectral sequence is established. The paper is neatly divided into three parts, the first based on Serre's methods, the second called ``The direct method'' based in Hochschild's manipulations with cochains, and the third establishes some interesting applications, in particular to the theory of simple algebras. 


The same techniques are applied without much pain to establish similar results for the cohomology theory of Lie algebras in \cite{kn:cla}.     
 
These constructions, were extremely important as they presented one of the first extensions of the very powerful computational tool of the spectral sequences introduced by Leray, to non--topological situations. Both sequences  were later encompassed by Grothendieck in his theory of derived functors in abelian categories; he proved the existence of a spectral sequence relating the derived functors of a composition, with the derived functors of its components. 

Around that time, Hochschild also started a line of work that he continued to pursue all along his later career, and that consisted in the generalization of the basic results in group and Lie algebra cohomology, to the situation of groups and algebras with additional structure where many of the standard techniques used in the discrete case do not apply. 

In \cite{kn:elg1} and \cite{kn:elg2} he tackles the problem of the classification of the extensions of Lie groups by constructing ``differentiable'' factor sets, a device that is not available in the situation of topological groups.

In the early 1950s started the interaction of Gerhard with the Bourbaki group, that began when he attended three of the group congresses even though he never was formally a member of the group.
In accordance with the Bourbaki files provided to the author by J.-P. Serre from Viviane le Dret, see \cite{kn:memorial}, 
Hochschild participated at the following instances: 
the congress at Pelvoux-le-Poet (June 25 to July 8, 1951), Foreign visitors: Hochschild and Borel, ``cobayes'': Cartier and Mirkil; the congress at Pelvoux-le-Poet (June 25 to July 8, 1952), Foreign visitors: Borel, De Rham and Hochschild; at the congress at Murols (August 17 to 31, 1954), Foreign visitors: Hochschild and Tate, ``Honorable foreign visitors'': Iyanaga and Yoshida; ``efficiency expert'' MacLane, ``cobaye'': Lang.

Hochschild mentioned later to the author, about his participation in the 1951 meeting, and the surprise that he felt when he saw the ``{\em very young boy with a boy scout look that went to pick me up at the station''}, later participating fully in the discussions of the group --of course he was referring to P. Cartier. 

It was in Urbana that Gerhard met his wife, Ruth Heinsheimer, who like Gerhard was born in Germany.  She and her mother escaped Germany in early 1939, first settling in Paris and then fleeing to a small village in the Pyrenees before sailing for New York from Lisbon in February of 1941.  Ruth graduated from Bryn Mawr College in 1947 and then enrolled in the graduate program in mathematics at the University of Illinois at Urbana -Champaing, where she obtained an M.A. degree in Mathematics in 1948.  When Gerhard arrived as an assistant professor, she was working under the direction of R. Baer. They married in July 1950 and their daughter Ann was born in 1955 in Urbana and their son Peter in 1957 in Princeton. Ruth passed away in El Cerrito, on June 2, 2005. 

In conversations with the author that took place much later, Gerhard recalled how much they both enjoyed very much those early years of their relationship at Urbana and frequently in his later years, mentioning nostalgically the loss of that close group of friends that included not only mathematicians but also some people in literature.


\subsection{The Berkeley period}

Gerhard remained at UIUC until September 1958, when he moved to Berkeley as a Professor of Mathematics, in the meantime Ruth had finished a Master's degree in mathematics and another in French literature. 

In the book {\em Mathematics at Berkeley, A History}, see \cite{kn:mab}, Calvin Moore presents a detailed and very insightful account of the process that starting from a minuscule department in a financially troubled private California college in the 1860s ended up producing one of the leading mathematical institutions in the world. 

In Chapter 13 of the book: {\em The Kelley years: 1957--1960}\,\footnote{See also C. Moore's contribution in \cite{kn:memorial}.} the author describes the strategy followed by John Kelley, then chairman of the department, to hire in clusters in order to develop the underrepresented areas: Algebra, Topology and Applied Mathematics.  The first cluster was to be in algebra where the department hired Hochschild and Maxwell Rosenlicht from Northwestern University. The deparment recommended the appointment of both of them as full professors effective July 1958 and made it known to each that the other was also being recommended. 


\medskip

Kelley knew Hochschild since they served together at Aberdeen, and Rosenlicht and Gerhard were acquainted since both were together at Harvard, the first as a graduate student and the second as a Pierce instructor. 

\medskip

The following year the same strategy was adopted in geometry with Spanier and Chern, and both accepted; Spanier arrived in 1959 and Chern in 1960. 


This went simultaneously with the hiring of some outstanding junior mathematicians: Glen Bredon, James Eels, Morris Hirch, Bertram Kostant, Emery Thomas and Steven Smale and altogether in {\em The Kelley years} the number of professorial staff increased from 19 to 44, with a very good balance between the fields. By 1960  the situation was much more balanced than in 1956 when the strength was basically in analysis (PDEs and functional analysis), computational number theory, and logic.


Hochschild had 26 students along his career and 22 of them graduated at Berkeley.  

Along his career, Gerhard's students thesis covered a wide thematic spectrum. Besides the expected subjects on homological algebra, Lie groups and Lie algebras, algebraic groups and Hopf algebras, one finds topics as Category theory, number theory and field theory, topological groups and groupoids, etc. 

His first PhD student was George Leger that finished in 1951 at the University of Illinois at Urbana--Champaign and wrote a thesis on {\em Cohomology theory for Lie algebras}, and his last student was Nazih Nahlus in Berkeley 1986 with a thesis entitled: {\em Lie theory of algebraic groups}. 

A constant theme in all comments of Hochschild's former students, is a deep appreciation of Gerhard's personality that went together with, but often transcended, his mathematical influence, mentioning his role as an advisor and later mentor and friend throughout their mathematical careers. 

Illustrative of the relationship he generated with his students is the comment of Ronald Macauley who wrote a thesis in 1955 at University of Illinois at Urbana on: {\em Analytic group kernels and Lie algebra kernels}.

In the acknowledgments he says: ``{\em The author takes this opportunity to express his most sincere gratitude to Professor G. P. Hochschild who offered and gave more aid and encouragement in the preparation of this thesis than could reasonably be expected of any advisor and friend. Indeed, Professor Hochschild's patience has withstood an arduous test''.}\footnote{We thank Prof. P. Bateman and the librarians of the University of Illinois for the information, see \cite{kn:memorial}.}

Hochschild frequently played the role of advisor of many junior members of the department, especially but not limited, to those in his areas of main interest. Even though he played a crucial academic role all along his tenure, he deeply disliked  the bureaucratic organization of academic institutions --for example he championed systematically and unsuccesfully, for the separation of sports and academics-- and for that reason he never accepted to take responsabilities in the administrative tasks of the University.

He served for a short time in a committee to select the new 
Instructors at Berkeley and the administration gave him guidance as to consider Affirmative Action policies. As soon as he received the information he told the Dean that he was unable to stay in the committee as he had personally seen the Nazis making a difference between ``jewish'' and ``aryan'' mathematics. He did not accept the Dean's arguments about the difference between ``negative'' racism as in the Nazi era and the policy they were implementing of ``positive'' racism. In a conversation with the author of this note, he mentioned that what was considered positive today, could easily become negative tomorrow and that the only safe attitude was to be blind to racial distinctions to the point of never registering any information about the issue, specially in official documents. He resigned from the committee.

While at Berkeley he was elected in 1979 to the National Academy of Sciences and to the American Academy of Arts and Sciences. In 1980 the American Mathematical Society awarded him the Steele prize for work of fundamental or lasting importance, citing in particular five of his papers published from 1945 to 1952 on homological algebra and its applications.

It was in the year 1980 that the author finished his PhD dissertation {\em Cohomology of comodules} working under his direction. In one of our frequent outings to drink coffee at the usual coffee shop in Hearst Ave., I asked him what he thought about his nomination for the Steele prize. His opinion was that to give him the prize was completely useless, he thought that the only utility a prize might have, was to help young people to get good jobs, and that for a senior mathematician like him, was completely irrelevant. 

During his tenure at Berkeley, Hochschild published 40 papers and all of his five books and along this very fruitful period of almost 30 years, his mathematical interest concentrated more and more in the theory of Lie groups, Lie algebras and algebraic groups, with special emphasis in the cohomological methods.  

His collaboration with Mostow blossomed with the stream of papers written on Representative functions (see in \cite{kn:memorial}, the articles by A. Magid and D. Mostow with a very rich description of the large ideas concerning Tannaka duality and representative functions; we also summarize the results in a later section).  

During the Berkeley period, and in a line of work with a strong homological flavor, he published the nowadays extremely cited paper \cite{kn:HKR}, where the so called Hochschild, Kostant, Rosenberg theorem appears. 

We briefly mention the importance that this result has had in the development of the basic ideas of ``non commutative geometry''. 

A convenient formulation of the HKR theorem is the following: for a smooth algebra $A$ over $\k$ there is a graded isomorphism $H_\bullet(A,A) \cong \Omega^\bullet_{A|k}$ where $H_\bullet(A,A)$ denotes the Hochschild homology with coefficients in $A$, and $\Omega^\bullet_{A|k}$ is the algebra of differential forms. The point here, taken over by Connes and others is that for a non commutative algebra $A$, one can think of $H_n(A,A)$ as the $n$--forms for the non commutative space described by $A$. Moreover in the ``geometric'' case of a smooth algebra, the Hochschild homology {\em is} the usual graded algebra of differential forms.  

Kostant was hired as an assistant professor at Berkeley in 1956 and left in 1961 to take a position at MIT. 

Both were interested in similar subjects in particular in the theory of algebraic groups that they studied together in the printed notes of the Chevalley Seminar 1956/1858, and as a result of their collaboration they wrote the paper \cite{kn:dflacag}, on the cohomology of homogeneous spaces for algebraic groups. In that paper they generalized the basic results on the rational cohomology of algebraic groups. This cohomology  theory was developed in \cite{kn:calg}. These two papers will be considered later in this note.

Kostant mentioned once to the author how much he agonized when he took the decision to leave Berkeley, mainly because the certainty of severely diminishing the relationship he had developed with Gerhard \footnote{See also the contribution by Kostant in \cite{kn:memorial}.}.

\subsection{Hochschild's mathematics at Berkeley}

We will comment at length in a later section on the contributions of Hochschild to the theory of algebraic groups and Hopf algebras. Now, we will concentrate briefly our attention on two papers on Lie theory, that concern specifically Ado's theorem and algebraic Lie algebras.

He produced a large opus in Lie theory, but many times along his life he went back in his efforts to some topics he had dealt with, at the beggining of his career. One of them has to do with Ado's theorem --and its global version as discussed for example in the contribution by Moskowitz in \cite{kn:memorial}.  

Ado's theorem, called sometimes Ado--Iwasawa's theorem,  
states that every finite-dimensional Lie algebra $\mathfrak g$  over a field $\k$ of zero characteristic, can be viewed as a Lie algebra of square matrices under the commutator bracket, i.e. $\mathfrak g$ has a faithful finite dimensional linear representation $\rho$ over $\k$. 
It was first proved in 1935 by Igor Dmitrievich Ado of Kazan State University, see \cite{kn:ado35}, and the restriction on the characteristic was removed later, by Iwasawa, Harish-Chandra and Hochschild, see \cite{kn:harish} and \cite{kn:iwasawa}.

After the original proof by Ado, valid for algebraically closed fields in characteristic zero and  that apppeared in 1935, E. Cartan published in 1938 an analytic proof in the case that the field was $\R$ or $\C$. The next published proof is due to Iwasawa 
in 1948, see \cite{kn:iwasawa},  and established the result in the case that the characteristic was different from zero. Soon afterwords, Harish--Chandra in \cite{kn:harish} --for characteristic zero-- published a proof that, besides being much more economical, was a little more precise than the original statement by Ado. 
In fact he proved that the faithful representation $\rho$ could be taken to have the additional property that every element of the maximal nilpotent ideal of $\mathfrak g$ is mapped into a nilpotent matrix on the representation space. 

It is interesting to note, that Harish--Chandra in his paper mentions in a footnote a proof by Hochschild: ``{\em Professor Chevalley has kindly informed me that Dr. Hochschild has succeded in constructing an algebraic proof which imitates Cartan's procedures''.}  The author of the present article, has never seen the original proof mentioned by Harish--Chandra (it seems that he did not publish it) but he mentioned once in a personal conversation, that his proof was the one adopted by Bourbaki in his treatise on Lie Groups and Lie algebras, Chapters 1--3.

In his 1966 paper, {\em An addition to Ado's theorem}, see \cite{kn:ado}, Hochschild presents a slight generalization of the classical theorem --as improved by Harish--Chandra. In the introduction he writes: ``{\em The main purpose of this note is to point out the following strenghtened (with respect to the nilpotency property) form of the theorem on the existence of a faithful finite--dimensional representation of a finite dimensional Lie algebra''.} Theorem 1 --that is the main statement of the paper-- reads as: {\em Let $L$ be a finite--dimensional Lie algebra over an arbitrary field, and let $\alpha$ denote the adjoint representation of $L$. There exists a finite dimensional representation $\rho$ of $L$ such that $\rho(x)$ is nilpotent for every element $x \in L$ for which $\alpha(x)$ is nilpotent.}\footnote{In the mentioned paper  Gerhard mentions Leonard Ross thesis written under his direction, {\em Cohomology of graded Lie algebras} for the suggestion that the ``{\em nilpotency property might be secured''.}}   

It is interesting to point out that Hochschild's proof --as almost all the proofs of these kind of results-- is not independent of the characteristic and in fact he produces two different proofs, one for the zero characteristic case and the other for positive characteristic.

It is obvious that Ado's theorem cannot be naively globalized, if $G$ is a Lie group and $\mathfrak g$ is its Lie algebra, there are serious obstructions to construct from the embedding  of $\mathfrak g$ into $\mathfrak {gl}_n$ for some $n$, an embedding of $G$ that linearizes to the original embedding of $\mathfrak g$. 

In the series of papers on representative functions written together with Mostow, the authors dealt with the globalization problem and proved the following neat result --see the contribution by Mostowski in \cite{kn:memorial}: 
Let $G$ be a connected real or complex Lie group. 
In the real case $G$ has a faithful finite dimensional smooth representation if and only if its radical and a Levi factor have such a representation. In the complex case $G$ has a faithful finite dimensional holomorphic representation if and only if the radical does (since a complex semisimple group always has a faithful representation).  If $G$ is either real, or complex with a faithfully represented Levi factor and a simply connected radical,  then $G$ has a faithful finite dimensional smooth (resp. holomorphic) which is unipotent on the nilradical --these results are also written in \cite{kn:book1}, Chapter XVIII. 

\medskip

Next, we discuss another important paper that appeared in 1971, called: {\em Notes on algebraic Lie algebras.} --see \cite{kn:algebraic}, where he studies a related problem. 

\medskip

In the introduction to the paper he writes as motivation : ``{\em  A Lie algebra is said to be algebraic, if it is isomorphic to the Lie algebra of an affine algebraic group. In view of the fact that entirely unrelated affine algebraic groups 
(typically, vector groups and toroidal groups) may have isomorphic Lie algebras, this notion of algebraic Lie algebra, calls for some clarification''}. If we start with an algebraic Lie algebra $\mathfrak g$, we call $G_{\mathfrak g}$ a connected affine algebraic group that has $\mathfrak g$ as Lie algebra--provided that it exists.
\medskip

A well known criterion for a Lie algebra to be algebraic is the so called Got\^o's theorem --published by M. Got\^o in 1948, \cite{kn:goto}-- that says that a finite dimensional Lie algebra $\mathfrak g$ over a field of characteristic zero is algebraic, if and only if its image under the adjoint representation is the Lie algebra of an algebraic subgroup of the group of automorphisms of $\mathfrak g$.   

\medskip

More precisely, consider $\mathfrak g$ and the representation $\operatorname {ad}: \mathfrak g \rightarrow \operatorname{D}(\mathfrak g) \subset \operatorname{End}(\mathfrak g)$, where $\operatorname {D}(\mathfrak g)$ denotes the Lie algebra of derivations of $\mathfrak g$. In this context $\operatorname {ad}(\mathfrak g) \subset \operatorname {D}(\mathfrak g)$ is the ideal of the inner derivations of $\mathfrak g$.  If we take the affine algebraic group $\operatorname{Aut}(\mathfrak g)$, then its associated Lie algebra is $\operatorname {D}(\mathfrak g)$, i.e. $\mathcal L (\operatorname{Aut}(\mathfrak g)) = \operatorname {D}(\mathfrak g)$.  In this situation Got\^o's theorem, asserts that it exists an algebraic subgroup $H \subset \operatorname{Aut}(\mathfrak g)$ with the property that $\operatorname {ad}(\mathfrak g)=\mathcal L(H)$ if and only if $\mathfrak g$ is algebraic.
In the case that the basic field is algebraically closed, 
$G_{\mathfrak g}$ can be taken to have unipotent center. 

Once that was unveiled the special unipotency property that one can enforce the center of the group $G_{\mathfrak g}$ to have, the road is paved for the following result, that is in fact the main theorem in \cite{kn:algebraic}. {\em Over an algebraically closed field of characteristic zero, there is exactly one isomorphism class of connected affine algebraic groups, with unipotent center, whose Lie algebras are isomorphic with a given algebraic Lie algebra.} --see \cite{kn:algebraic}, page 10. 

The situation is simpler in the case of unipotent groups and nilpotent Lie algebras. With the methods developed in the papers on representative functiones specially in 
{\em Algebraic groups and Hopf algebras} ,\cite{kn:agha}, it can be proved that  the category of  unipotent affine algebraic groups over a field of characteristic zero is equivalent to the category of finite dimensional nilpotent Lie algebras.

Using Got\^o's theorem and the above result on nilpotent Lie algebras, Hochschild establishes his main result by first decomposing the group $G_{\mathfrak g}$, given by Got\^o, into the semidirect product of a linearly reductive group and the unipotent radical and treating both cases separatedly.  

He also mentions that the result on isomorphisms cannot be strengthened to morphisms and exhibits some examples that show that a Lie algebra homomorfism of algebraic Lie algebras is not always the differential of a homomorphism of the affine algebraic groups constructed from the Lie algebras, even if the groups have unipotent centers.

The simplest example is the following: let $\k$ be an algebraically closed field of characteristic zero and $G=\k_a \rtimes \,\k_m$ (additive and multiplicative groups of the field respectively), with product $(a,u)(b,v)=(a+ub,uv)$, then $\mathcal L(G)= \k x + \k y$ with $[x,y]=y$ and $\mathcal L(\k_a)= \k y$. The Lie algebra morphism $\rho: \mathcal L(G) \rightarrow \mathcal L(\k_a)$ that sends $\rho(x)=y, \rho(y)=0$ is surjective but it is not the differential of a homomorphism of algebraic groups because there is no morphism of $\k_m$ onto $\k_a$.  

Moreover, the above result proves that the natural map of the automorphism group of such an algebraic group --with unipotent center-- into the automorphism group of the Lie algebra, is an isomorphism. This implies that the group of algebraic automorphisms of such kind of groups, is also an affine algebraic group. This situation is extremely rare, the only other example of such a family --as Hochschild himself proved in another paper, is the family of groups with the property that the dimension of the centers of $G/G_u$ is zero or one --see Chapter XV of Hochschild's book {\em Basic theory of algebraic groups and Lie algebras}, \cite{kn:book1}.
 
\subsection{Hochschild's books}

\begin{enumerate}

\item {\em Perspectives of elementary mathematics}. New York-Heidelberg-Berlin: Springer-Verlag.(1983).

\item {\em Basic theory of algebraic groups and Lie algebras}. Graduate Texts in Mathematics, 75. New York-Heidelberg-Berlin: Springer-Verlag. (1981).

\item {\em Introduction to affine algebraic groups}. San Francisco-Cambridge-London-Amsterdam: Holden-Day, Inc. (1971).

\item {\em A second introduction to analytic geometry}. San Francisco-Cambridge-London-Amsterdam: Holden-Day, Inc. (1968).

\item {\em The structure of Lie groups}. Holden-Day Series in Mathematics. San Francisco-London-Amsterdam: Holden-Day, Inc. (1965).
\end{enumerate}

Gerhard wrote the five books mentioned above along his career all of them while at Berkeley. 

Three of them: 
{\em Basic theory of algebraic groups and Lie algebras}, 1981; {\em Introduction to affine algebraic groups}, 1971 and {\em  The structure of Lie groups}, 1965, are addressed primarily to graduate students, and the other two have  more global educational goals. 

The purpose of the book: {\em Perspectives of elementary mathematics} becomes clear in the handwritten dedication that Gerhard stamped in this author's copy: ``{\em You may view this as my anti--education manifesto}''. In the Preface, Gerhard writes: {\em The general aim of what follows is to present basic mathematical concepts and techniques in familiar contexts in such a way as to illuminate the nature of mathematics as an art...In order to avoid burying the essentials under routine technicalities, a style has been adopted that relies on the reader's active involvement somewhat more than it is customary..''}

A novelty of the presentation of the material, is that at the end of each chapter there is a project concerned with the creation of computer programs --without the need to do any formal programming. For example, at the end of Chapter IX entitled: {\em The sphere in 3--space}, where in particular he exhibits the identification (up to a sign) of the group of unit quaternions with the group of rotations of $\R^3$, he suggests the following project: ``{\em Design computer routines implementing quaternion arithmetic. There should be four functions  of quaternions with quaternions as values: sum, negatives, product, reciprocal. Note that such a facility contains complex number arithmetic.}'' 

\medskip

The book: {\em A second introduction to analytic geometry}, has a somewhat similar purpose. It was dedicated to his son Peter when he was an undergraduate, and in the preface Gerhard wrote: ``{\em What follows is an examination of the basic geometrical features of Euclidian three--space from the point of view of rigorous mathematics. This requires that even the simplest visually obrious facts concerning the relations between points, lines, and planes in space be rigorously deduced from the properties of the system of real numbers\ldots [and] indeed, our program here necessarily involves algebra and analysis as much as it involves geometry''.}

\medskip
Next we take a quick glance at the  other three books, that have been designed as texts for graduate students.

The book, {\em The structure of Lie groups} received a very serious attention from reviewers and specialists as some of its reviews clearly show.

\medskip

The titles of the Chapters indicate by themselves the extremely ambituous scope of the book: 1. Topological groups, 2. Compact groups, 3. Elementary structure theory, 4. Coverings, 5. Power series maps, 6. Analytic manifolds, 7. Analytic subgroups and their Lie algebras, 8. Closed subgroups of Lie groups, 9. Automorphisms groups and Semidirect products, 10. The Campbell--Hausdorff formula, 11. Elementary theory of Lie algebras, 12. Simply connected analytic groups, 13. Compact analytic groups, 14. Cartan subalgebras, 15. Compact subgroups of Lie groups, 16. Centers of analytic groups and closures of analytic subgroups, 17. Complex analytic groups, 18. Faithful representations. 

\medskip

The above list of wide ranging but rather minimalist titles, does not explicitly illustrates one of the great strengths of the book. Actually, besides providing unified proofs for many important and then recent results that were scattered in seminar notes or specialized articles, not few of 
new and not previously published results are included. For example: In  Chapter 15 the author presents a deep generalization of a theorem due to Iwasawa in order to obtain the basic results concerning the maximal compact subgroups of Lie groups.  Moreover, in Chapter 18 his results with Mostow on the globalization of Ado's theorem mentioned before, concerning the existence of finite  dimensional representations of analytic groups are exposed for the first time in book form. 

\medskip 

K. H. Hofmann, in his review for Zentralblatt MATH, 1966, Zbl 0131.02702, comments on the need for a book on the subject that: ``{\em presents the foundations of [global] Lie theory in a language that takes into account the recent developments'',} and gives his opinion that: ``{\em The book by G. Hochschild is an admirable response to that need''}. He commends the fact that:  ``{\em The book is completely self--contained insofar as not results are used whose proofs have to be looked up elswhere. The prerequisites are held at a minimum in order not to discourage interested students\ldots The whole architecture of the book is very clear and systematic. It is amazing to see the author present such a large body of information on 226 pages''.}

\medskip

Not everybody found satisfactory the minimalist style of the book, and some authors described it as ``relentless'', even though other reviewers did not share at all the negative opinion.  Prof. F. Hirzebruch writes: ``{\em It is amazing the extent to which the author achieved his goal to enable a self contained reading to someone who only knows the basics of multilinear algebra, group theory, set theoretical topology and calculus''.}

\medskip

Hochschild wrote two books on the theory of algebraic groups, the first: {\em Introduction to affine algebraic groups} in 1971 and the second: {\em Basic theory of algebraic groups and Lie algebras} in 1983. 

In the period of a few years, and starting with A. Borel seminal book: {\em Linear algebraic groups}, that appeared in 1969 and not counting different sets of seminar notes, the mathematical community saw the appearence of a flood of books on the theory of algebraic groups, that reflected the deep interest of the mathematical community on the subject.
 
The following might be a not too incomplete list: M. Demazure and P. Gabriel: {\em Groupes Alg\'ebriques}, 1970; E. Kolchin: {\em Differential algebra and algebraic groups}, 1973; J. Humphreys : {\em Linear algebraic groups},  1975; W. Waterhouse: {\em Introduction to affine group schemes}, 1979; D.G. Northcott: {\em Affine sets and affine groups}, 1980, T. Springer: {\em Linear algebraic groups}, 1981.  All of these monographs have as predecessors three basic references of Chevalley: {\em Th\'eorie des groupes de Lie II, Groupes Alg\'ebriques}, 1951; {\em Th\'eorie des groupes de Lie III, Groupes Alg\'ebriques}; 1954, {\em Classification des groupes de Lie alg\'ebriques}, 2 vols, Notes polycopi\'ees, 1956/58 (see also the edition of Chevalley's collected works, \cite{kn:checartier}).  Of these three works of Chevalley, the last has a very different flavor from the first two: the methods have changed radically as algebraic geometry marched into the subject after the work of Borel and Chevalley himself, and the algebraico--geometric emphasis --as contraposed with the characteristic zero, Lie algebra--Lie group inflexion-- is very explicit in the majority of the books of the above list.   

Hochschild's two books cover different material than the others, and there is less emphasis in the classification theorems of semisimple algebraic groups (that are the main topic of Chevalley's seminar and justly considered a major --and also slightly suprising-- achievement of the theory) and there is more emphasis in his own contributions to the subject --many of them jointly with Mostow.

In that respect we mention the following subjects (without distinguishing in which of the two books they are included): Mostow's theorem about the  decomposition of an affine group in characteristic zero as a semidirect product of its radical and a linearly reductive group; the result mentioned above on algebraic Lie algebras and its consequences concerning necessary and sufficient conditions for the automorphism groups of an affine algebraic groups to be affine; the theory of observable subgroups as developed jointly with Mostow and Białynicki--Birula; the theorems about the categorical equivalence of the nilpotent Lie algebras with the unipotent groups, etc. 

But what more sharply distinguishes Hochschild's books from the others, are the following two features --one of technical and the other of stylistic character, the first is the systematic use of Hopf algebra theory in order to get control of the basic features of the affine theory, the second the strong determination to keep the monographs as much self contained as possible. 

The second feature is not surprising and can be considered as his usual style of writing mathematics--see above. For example, in the introduction to the 1971 book, he writes: ``{\em In order to keep this book self--contained, we have limited the material so that only elementary affine algebraic geometry comes into play, and all the [needed] results in this area\ldots are established in Section 1\ldots no special knowledge is presupposed, although it is assumed that the reader will be able to assimilate the daily diet of the working algebrist in unpremasticated form''.} The same can be said about the second book. For that reason the books have a manageable size, and can be used for standard courses in the subject. 

Concerning the first issue, it is worth mentioning that in the opinion of some integrants of the mathematical community the Hopf algebra viewpoint that Hochschild adopts in a rather punctilious manner ``{\em seems to detract from the geometric content of the elementary results\ldots}'' c.f. \cite{kn:humphreview}. 

The author of this article does not intend to enter in the rather epistemological discussion that would be necessary in order to try to understand expressions like: ``geometric content'' or ``geometric methods'', specially when they are used in opposition to ``algebraic methods'', but we find useful to point at a few aspects of the theory where the Hopf algebra viewpoint adopted by Hochschild, simplifies considerably the treatment of the subject\footnote{To be fair, we must mention that the Hopf algebra viewpoing has been used systematically in the author's book (written together with A. Rittatore), ``Actions and invariants of algebraic groups'', 2005, see \cite{kn:libro}}.  

To deal with quotients in the category of algebraic varieties is more complicated than in other geometrical environments because of the scarcity of functions we have at our disposal. 

In particular, if $G$ is an affine algebraic group and $K$ a closed normal subgroup, the quotient group has a natural structure of affine algebraic group. This is a classical result --probably due to Chow in the 50's-- and the known proofs are in general non intrinsec and rather indirect. 

In his 1971  book, Hochschild presents this result in the following not hard to prove intrinsec form. He proves that if $A$ is an affine Hopf algebra and $B$ a sub Hopf algebra, then $B$ is also affine. Taking $A=\k[G]$ and $B=\k[G]^K$ where $K$ is a normal subgroup of $G$, we obtain the quotient by dualizing $B$.

Moreover, the systematic use of the coproduct $\Delta: A \rightarrow A \otimes A$, induced on $A=\k[G]$ by the product on the affine group $G$, simplifies some formulas and definitions, trivializing the operational aspects. 
 
For example if $\tau,\sigma: A \rightarrow \k$ are two tangent vectors at the identity, its Lie bracket is simply the commutator associated to the convolution product: $[\tau,\sigma]=\tau \star \sigma - \sigma \star \tau = (\tau \otimes \sigma - \sigma \otimes \tau) \Delta$, and if  $\sigma: A \rightarrow \k$ is a tangent vector at the identity as above, the associated vector field --invariant derivation-- can be expressed simply as: $(\sigma \otimes \operatorname{id})\Delta: A \rightarrow A$.

It is this author's opinion, that Hopf algebra theory plays in this area the role that differential calculus plays in the theory of Lie groups, and plays it better as it is more algorithmic. Hence, one can view the systematic introduction of the Hopf algebra structure as one more step in the process of ``{\em algebraiz[ing]\ldots to death}'' the theory of algebraic groups (Chevalley's ditto) --this expression was taken from the very interesting monograph by A. Borel, see \cite{kn:borelhistory}, Chapter 7, page 147.

\subsection{Retirement and photography}

Hochschild retired on July 1, 1982, teaching part-time until retiring fully on July 1, 1985\footnote{See the article by C. Moore in \cite{kn:memorial}.}.

This happened in accordance to the general retirement policy of the university that required mandatory retirement of tenured faculty on July 1 following their $67^{\text{th}}$ birthday. He always resented this policy that soon changed back and forth with the waves of US national politics.   

After retirement he dedicated much time to two of his all life passions, cinema and photography.

He particularly enjoyed the not mainstream movies that were shown in ``cin\'emath\`eques'' and alternative movie theaters in the Bay Area, where he used to go --by himself or with some of his friends, ex--students or colleagues, and enjoy the movies, sitting in a corner far away from the crowds. I remember an all night stand where we watched a whole series of the Clint Eastwood's spaghetti westerns. He mentioned then, that in other of these marathons, he had seen the same series accompanied by Moss Sweedler. 

A more important area of his attention was photography. He was keen on the activity since his childhood in Berlin, and later in life while at Berkeley he took long trips in the west --including Canada and Alaska-- to take pictures. He was already taking regular photography trips in the 60s and after retirement he became more and more involved in landscape photography, pursuing this hobby with tenacity. He would go off on long trips by himself --at the end accompanied by a cell phone and once he told this author, that in the places he went he hardly ever had coverage-- with his gear, driving thousands of miles to find the right light and framing. Each picture took him at least one day and he amassed an incredible amount of negatives that he himself developed at the lab he had built at the basement of his house. His pictures were almost always black and white and he constructed at his house a device, that he installed in the living room wall close to the piano and to the his go board, where he displayed a changing exhibition for visitors. 



He died peacefully in El Cerrito, accompanied by his daughter Ann at his bedside on July 8, 2010, at the house where he lived since he moved to the Bay Area in 1958.  
\section{Hochschild's work on algebraic groups and Hopf algebras}

Hochschild work on algebraic groups is vast. He wrote dozens of papers on the subject and it won't be possible to present a detalied analysis of his whole opus. We will concentrate upon three aspects of his work: we consider his contributions to representative functions of Lie and analytic groups; his work on representation and cohomology of groups with additional geometric structure that includes linear algebraic groups; and on his one paper in invariant theory of unipotent groups. Other relevant aspects of his production like the study of automorphism groups of affine groups and profinite groups will be omitted.
\subsection{From Lie groups to algebraic groups}
Hochschild involvement with the general theory of affine algebraic groups\footnote{In the note by Kostant that will appear  in \cite{kn:memorial} the author recalls that at the late fifties when both were together at Berkeley, he and Gerhard used to meet to expose to each other the material in the Chevalley's seminars: on the foundations of algebraic geometry \cite{kn:chefund} and on the classification of semisimple algebraic  groups \cite{kn:checlassif}.} is closedly related with Chevalley's work on the subject, for that reason we start with a few words about Chevalley's opus in that period. 

In Chapters VI and VII of the monograph {\em Essays in the History of Lie Groups and Algebraic Groups}, by A. Borel --see \cite{kn:borelhistory}-- the author recounts that the starting point by Chevalley was in his 1943 paper: {\em A new kind of relationship between matrices}. There he generalized to general fields Maurer's results on (the later called), algebraic Lie algebras, proved over $\C$. For that, he defined the notion of ``replica'' of a matrix and proved that Maurer's result can be expressed as follows: {\em If the Lie algebra of a complex linear algebraic group contains a matrix, it also contains all its replicas}. Later in a joint paper with Tuan, 1946, they proved the converse. All this was later incorporated in \cite{kn:che2} and \cite{kn:che3}.

In Chevalley's classic 1951 book: {\em Theory of Lie groups}, Chapter VI: {\em Compact groups} --see \cite{kn:che1}, he deals with Tannaka duality theory --see \cite{kn:tannaka}-- and the approach he presented was followed --and vastly generalized-- in the series of papers on Representative functions by Hochschild and Mostow, see \cite{kn:repre1}, \cite{kn:repre2}, \cite{kn:repre3}. The authors of the series mention also other relevant work besides Chevalley's that they used: one is the paper by Harish--Chandra in 1950, see \cite{kn:harish2}, and the other is the 1956 paper by Cartier, see \cite{kn:cartier}\footnote{See also the personal reference to this topic that appears in the interview to P. Cartier in \cite{kn:memorial}.}. 

We start by explaining the basic ideas of Chevalley--Hochschild--Mostow approach to Tannaka's theory in modern guise --we have used heavily the very illustrative survey by Joyal and Street: {\em An introduction to Tannaka duality and quantum groups, \cite{kn:js}}.


Assume that $G$ is a compact topological group and call ${}_G\mathcal M$ the monoidal rigid category of its finite dimensional complex representations. For a finite dimensional $G$--module $V$ we call $\rho_V:G \rightarrow \operatorname{GL}(V)$ the morphism associated to $V$.  

This monoidal category is fibered over $\mathcal V= \mathcal V_\C$ (the category of finite dimensional complex spaces) via the forgetful functor $U: {}_G\mathcal M \rightarrow \mathcal V$.

One can take $\mathcal E(U)$ the set of natural transformations of $U$ that when equipped with the vertical composition and the natural addition, becomes an algebra. Moreover, it can be endowed with a natural topology that makes it a topological algebra: we take the coarsest topology that makes all the proyections $\mathcal E(U) \rightarrow \operatorname{End}(V)$ continuous for all $V \in {}_G\mathcal M$.  Recall that $\mathcal E(U)$ has a conjugation  operation given as follows. For any  natural transformation $u$, we define the natural transformation $\bar{u}_V= (u_{V_c})_c$ where $(-)_c: {}_G\mathcal M \rightarrow {}_G\mathcal M$ is the standard conjugation functor; i.e. the functor that associates to a vector space the same abelian group but the multiplication by scalars is defined by composing conjugation of the scalar and then applying the original multiplication.  

Next we restrict the natural transformations to the monoidal ones, i.e., we assume that $u_{V \otimes W}= u_V \otimes u_W$ and $u_\C=\operatorname{id}_\C$. 

Clearly, if we take only the natural {\em monoidal} transformations, we have a closed submonoid --called $\mathcal E^\times(U)$-- of $\mathcal E(U)$ that is in fact a topological monoid with the topology induced by the one of $\mathcal E(U)$.  

It is not hard to prove that in the case that $G$ is a group the monoidal transformations are invertible. Indeed if $u \in \mathcal E^\times(U)$ and we define $u' \in \mathcal E(U)$ as $u'_V= (u_{V^\vee})^\vee: V \rightarrow V$, where $V^\vee$ is the (contra-gradient) representation given by the rigidity property of the category ${}_G\mathcal M$, then, $u'$ is also monoidal and the inverse of $u \in \mathcal E(U)$.  

\begin{defi} In the context above, given the topological group $G$ we write $A(G) = \mathcal E^\times(U)  \subset \mathcal E(U)$. Thus, $A(G)$ is the group of all monoidal natural transformations of the forgetful functor on ${}_G\mathcal M$ with the topology induced by the one of $\mathcal E(U)$. We define the closed subgroup  $B(G)= \{u \in A(G): u=\bar{u}\}$. 
We avoid the name of Tannaka group for
$B(G)$ that is used in \cite{kn:js} as it is not standard. The group $A(G)$ has been called by Lubotzki the Hochschild--Mostow group. Hochschild and Mostow called it the group of proper automorphisms of $\mathcal R(G)$, we explain the reason for that nomenclature later.
\end{defi}

Clearly, if $x \in G$ the left translation by $x$ on any representation of $G$,  defines a monoidal natural transformation. We call this function $\pi : G \rightarrow B(G)$. 

One has the following theorem.

\begin{theo} The morphism $\pi: G \rightarrow B(G)$ is a continuous homomorfism of topological groups and if $G$ is compact so is $B(G)$. 
\end{theo}

The compactness of $B(G)$ follows from the fact that we can use a Haar measure to endow all the representations with an invariant inner product and then deduce from the naturality, that all the components $u_V$ are unitary operators with respect to that inner product.

Next we define the general concept of representative function. Assume, that $G$ is an abstract group and take $\C$ as the base field.    

\begin{defi} A function $f:G \rightarrow \C$ is said to be a representative function, if for all $x \in G$, the set of left translates $\{x\cdot f: G \rightarrow \C: \forall x \in G\}$ spans a finite dimensional space in the vector space of all functions from $G$ into $\C$. Recall that in the above context, one defines the left and right translations of a given function by the formulas: $(x \cdot f)(y)=f(yx)$, $(f \cdot x)(y)=f(xy)$.
\end{defi}

One easily sees that if the finiteness condition holds for the left translates, it also holds for the right translates and even for the two sided translates. It is clear that the set of all representative functions of $G$ forms a commutative $\C$--subalgebra of the algebra of all functions, and that it is closed by conjugation. It will be denoted as $\mathcal R(G)$.

Even if, it was never explicitly displayed in the early references, it is clear --and nowadays well known-- that $\mathcal R(G)$ has a natural Hopf algebra structure over the base field\footnote{In the interview to P. Cartier that appears in \cite{kn:memorial}, he mentions that it was in their 1958 meeting at Urbana, that he probably suggested the importance of Hopf algebra theory in relation with the representative functions, introducing Gerhard to that notion.}. We spell out this viewpoint in order to show how deeply embedded the Hopf algebra viewpoint was in Hochschild's early work, even if it was only years later that became more explicit.

The basic finiteness property of a representative function, guarantees that for $f \in \mathcal  R(G), x \in G$, one can find $f_i,g_i \in \mathcal R(G)$ such that $x\cdot f=\sum f_i(x)g_i$, and the morphism $\Delta: \mathcal  R(G) \rightarrow \mathcal  R(G) \otimes \mathcal  R(G)$ defined as $\Delta(f)=\sum f_i \otimes g_i$ is a coproduct in $\mathcal  R(G)$, compatible with the product. We will adopt Sweedler's notation and write $\Delta(f)=\sum f_1 \otimes f_2$. It is easy to show that the map $\mathcal S(f)(x)=f(x^{-1})$ defines an antipode, and that the evaluation at the identity of the group defines a counit. Similarly, to an arbitrary representation $V$ of $G$, we associate an $\mathcal  R(G)$ comodule structure $\chi_V: V \rightarrow V \otimes \mathcal  R(G)$ by the formula $\chi_V(v)=\sum v_0 \otimes v_1$, if and only if for all $x \in G, v \in V$, $x\cdot v=\sum v_0 v_1(x)$. 

In case that the group $G$ has additional structure, it is customary to ask the representative functions to preserve that additional structure, for example for a topological group, one asks the functions to be continuous, i.e. $\mathcal R(G) \subset \mathcal C(G,\C)$, etc.  

We concentrate from now on in the case of a toplogical group $G$, and in this situation  $\mathcal R(G)$ is a closed by conjugation subalgebra of $\mathcal C(G,\C)$ that has a natural Hopf algebra structure. 

The famous Peter--Weyl theorem, asserts that if $G$ is a compact topological group then $\mathcal R(G) \subset \mathcal C(G,\C)$ is a dense subalgebra with respect to the uniform topology.

It is usual to consider for an arbitrary group $G$ the property of having sufficiently many representations. 
The group $G$ has sufficiently many representations if for all $1 \neq x \in G$, there exists a finite dimensional $G$--module $V$, such that $\rho_V(x) \neq \operatorname{id}$. In other words, in $V$ there is an element $v \in V$ such that $x\cdot v \neq v$. It is easy to prove that $G$ has sufficiently many representations if and only if $\mathcal R(G)$ separates the points of $G$. 

Using Peter--Weyl theorem, it can be shown that a compact group has sufficiently many representations.

\begin{defi}
Assume that $\mathcal X \subset {}_G\mathcal M$ is a family of representations. We say that $\mathcal X$ is closed if it is closed under: isomorphisms, subrepresentations, direct sums, tensor products, conjugation, and $\C \in \mathcal X$. 
\end{defi}

For any family of objects $\mathcal Y \subset {}_G\mathcal M$, we call $\mathcal R(\mathcal Y)$ the family of representative functions corresponding to the elements of $\mathcal Y$. In particular $\mathcal R({}_G\mathcal M)=\mathcal R(G)$.

The proof of the assertion that follows is not short and requires the use of the orthogonality relations. It is written down in \cite{kn:che1}, Chapter VI.

\begin{theo} \label{theo:class} Assume that $G$ is compact and that $\mathcal X$ is a closed family of representations with the following property, for every $x \neq 1 \in G$, there exists $V \in \mathcal X$ and $v \in V$ such that $x\cdot v \neq v$, then $\mathcal R(\mathcal X)=\mathcal R(G)$. 
\end{theo}

We have all the machinery ready to illustrate the proof of the theorem of Tannaka--Krein following Chevalley's methods.
\begin{theo}
Assume that $G$ is a compact topological group, then $\pi: G \rightarrow B(G)$ is an isomorphism.
\end{theo} 

We proceed as follows: there is a natural restriction functor $\pi^*: {}_{B(G)}\mathcal M \rightarrow {}_G\mathcal M$ and a natural extension functor $e: {}_G\mathcal M\rightarrow {}_{B(G)}\mathcal M$, that given $V \in {}_G\mathcal M$ sends it into the representation of $B(G)$, that sends $u \in B(G)$ into $u_V:V \rightarrow V$. One easily shows that $e$ preserves sums, tensor products, conjugate representations, irreducibility and neutral element. This guarantees that the image of $e$ is closed as a family of representations of $B(G)$. Now, if $1 \neq u \in B(G)$, for some representation $V$ of $G$, $u_V \neq \operatorname{id}_V$ and then there is an element $v \in V$ such that $u \cdot v =u_V(v) \neq v$. Hence $\mathcal X$ has all the representations of $B(G)$ and then $e$ and $\pi^*$ are equivalences. It follows easily from the above that the morphism $\pi^*: \mathcal R(B(G)) \rightarrow \mathcal R(G)$ is an algebra isomorphism.  

To finish the proof of Tannaka--Krein we observe first that $\pi$ is injective. It follows from the Theorem of Peter--Weyl that for any $1 \neq x \in G$, there is a representation $V$of $G$ and an element $v \in V$ such that $x \cdot v \neq v$. That means that $\pi(x) \neq 1$. The surjectivity of $\pi$ is slightly more laborious and one has to use in an appropriate manner, the orthogonality relations and trace formula. 

\medskip

Introducing the ideas of Tannaka reconstruction --that we do not explain here-- one can go one step further and prove that the real spectrum of $\mathcal R(G)$ is isomorphic to $G$, when $G$ is compact. We call $\operatorname{Spect}_\R(\mathcal R(G)) \subset \mathcal R(G)^*$ the set of elements of the dual that are multiplicative and commute with conjugation. 

We start by defining the {\em generalized Fourier transform} , that is a morphism $\mathcal F: \mathcal R(G)^* \rightarrow \mathcal E(U)$, given for $\alpha \in \mathcal R(G)^*$ as: $\mathcal F(\alpha)_V : V \rightarrow V$, $\mathcal F(\alpha)_V= (\operatorname{id} \otimes \alpha)\chi_V$.  It is clear that if we call $\operatorname{ev}: G \rightarrow \mathcal R(G)^*$ the natural evaluation morphism, one has that the triangle below commutes.  

\medskip
\begin{equation}
\label{eqn:basic}
\xymatrix{& \mathcal R(G)^* \ar[dr]^{\mathcal F} & \\
G \ar[ur]^{\operatorname{ev}} \ar[rr]^{\pi}& & \mathcal E(U).}
\end{equation}

It can be proved --by describing the inverse of $\mathcal F$ a map that is sometimes called the Fourier cotransform-- that: $\mathcal F(\operatorname{Spect}_\R(\mathcal R(G)) = B(G)$, and that with  adequate topologies $\mathcal F$  is in fact an isomorphism of topological groups.

In that case we have a diagram as below, where all the arrows are isomorphisms:

\medskip
\begin{equation}
\label{eqn:basic2}
\xymatrix{& \operatorname{Spect}_\R(\mathcal R(G)) \ar[dr]^{\mathcal F} & \\
G \ar[ur]^{\operatorname{ev}} \ar[rr]^{\pi}& & B(G).}
\end{equation}

More information can be obtained if the group $G$ can be faithfully represented. Next lemma shows that the compact groups that can be faithfully represented verify that $\mathcal R(G)$ is finitely generated, in fact it is an affine algebra.

\begin{lema} \label{lema:finitegeneration} Assume that $G$ is compact and admits a faithful finite dimensional representation $V$. Then $
\mathcal R(G)$ is generated by the matrix coefficients of $V$ and the inverse of the determinant function on $V$. 
\end{lema} 

Call $\mathcal A$ the subalgebra of $\mathcal R(G)$ generated by the matrix coefficients of $V$ and the inverse of the determinant and call $\mathcal X$ the collection of all representations whose coefficients belong to $\mathcal A$. One proves that $\mathcal X$ is a closed family and as $V$ is faithful we can apply Theorem \ref{theo:class} and conclude that $\mathcal X = {}_G\mathcal M$ and then that $\mathcal A = \mathcal R(G)$. 

\begin{theo} \label{theo:surprising} If $G$ is a compact Lie group, then $G$ is real algebraic.
\end{theo}

This result follows by putting together the facts proved above: the well known fact that a compact Lie group admits a faithful representation and Lemma \ref{lema:finitegeneration} guarantee the finite generation of $\mathcal R(G)$, then by the fact that the map $\operatorname{ev}$ appearing in diagram \eqref{eqn:basic2} is an isomorphism, we conclude that the $G$ is indeed real algebraic. 

\medskip

The above slightly surprising Theorem \ref{theo:surprising}, was the ultimate reason that led the group of mathematicians that were working around Chevalley in the theory of Lie groups, to start switching their attention to the study of algebraic groups.

As we mentioned before, at the end of the 1940s, Chevalley himself started to work in the foundational basis of the theory of linear groups\footnote{Some of the other mathematicians that were working in this theory at that time were: Barsotti, Borel, Cartier, Chow, Dieudonn\'e, Freudental, Got\^o, Grothendieck, Harish--Chandra, Kolchin, Lang, Matsushima, Rosenlicht, Serre, Springer, Tits, Weil, etc.}. In the second and third books of the series {\em Theory of Lie groups, I,II,III}, see \cite{kn:che1,kn:che2,kn:che3}, he concentrated more fully on algebraic groups. In the second, he dealt with algebraic groups over arbitrary fields, defining them as subgroups of a general linear group. Seen from today's  viewpoint, the foundational principles and language he adopts does not seem to be the more efficient and effective, for example, it is not easy to define adequately the notion of homogeneous space or even of quotient group. Moreover, the main results are established in characteristic zero, and this restriction is forced because of his systematic use of the formal exponential, a device that only makes sense in that case. But, in any case he manages to develop many the essential foundational ingredients: for example at the end of the book he presents a proof of the existence and uniqueness of the multiplicative Jordan decomposition of the elements of the group (called sometimes today the Jordan--Chevalley decomposition).  The third book of the series, \cite{kn:che3}, is a mixture of Lie and algebraic group theory. He concentrates in the study of Lie algebras and in particular establishes the main results on the conjugacy of Cartan subalgebras and defines the concept of Cartan subgroup. The proof that a compact Lie group is real algebraic that we mentioned above, and that is hinted in the first book, is explicitly developed here. We extracted the above considerations from \cite{kn:borelhistory}, where a  thorough and authoritive description of the process mentioned above is presented.
 
\medskip

Let us briefly describe the first three of Hochschild's 
papers with Mostow on representative functions, see \cite{kn:repre1,kn:repre2,kn:repre3}.

The first of these papers was written while both were at the Princeton Institute during the year 1956/57, and in the second and third Gerhard gives as address, University of Illinois and University of California, respectively.

The authors build upon the foundations laid out by Tannaka in \cite{kn:tannaka} and Chevalley in \cite{kn:che1}, but with certain important modifications, that allowed them to work with much more flexibility. 

First of all, while in the original constructions Chevalley worked with the ``characters'' of the algebra $\mathcal R(G)$, i.e. with the maximal spectrum, Hochschild and Mostow choose to work with the proper automorphisms, i.e. the algebra morphisms of $\mathcal R(G)$ that commute with the right translations. Of course, both families of maps are in bijective correspondence but, the algebraic structures on the set of automorphisms are more natural. In modern lenguage and given the comultiplication $\Delta: \mathcal R(G) \rightarrow \mathcal R(G) \otimes \mathcal R(G)$, and the counit $\varepsilon$, to a character $\sigma: \mathcal R(G)\rightarrow \C$ one associates the proper automorphism   
$X_\sigma= ( \operatorname {id} \otimes \sigma)\Delta:\mathcal R(G) \rightarrow \mathcal R(G) \otimes \mathcal R(G) \rightarrow \mathcal R(G)$. Conversely given a proper automorphism $X:\mathcal R(G) \rightarrow \mathcal R(G)$ it has an associate character $\sigma_X = \varepsilon X: \mathcal R(G) \rightarrow \C$. We only need an algorithmic manipulation to prove that these correspondences are inverses of each other.  
 
In this situation, and adapting the notations already introduced, we call now --as the authors do-- $A(G)$ the group of proper automorphisms, and $B(G)$ the group of proper automorphisms that commute with complex conjugation. The fact that these groups correspond to the ones previously defined via natural transformations was explained above, and the identification is given by the Fourier transform. In this context, it is clear that the morphism $\pi: G \rightarrow A(G)$ considered before is simply the left translation on $\mathcal R(G)$.   

\medskip

Another crucial ingredient that \cite{kn:repre1} brings to the theory, is the clarification of the role of the base field, and this is achieved via the concept of {\em universal complexification.} If $G$ is an arbitrary real Lie group, they define a complex Lie group $G^+$, called the universal complexification of $G$ that comes equipped with a canonical morphism $\iota: G \rightarrow G^+$ and its basic property being that there is a bijective correspondence between the complex representations of $G$ and the complex representations of $G^+$. If $\rho$ is a complex representation of $G$ and we call $\rho^+$ the associated complex representation of $G^+$, then $\rho = \rho^+ \iota$. In particular it is clear that $\pi: G \rightarrow A(G)$ extends to $\pi^+:G^+ \rightarrow A(G)$. 

\medskip

One needs at this point --for deep reasons that have to do with the fact that contrary to the complex case, in the real situation an irreducible real algebraic group need not be connected-- to impose a hypothesis of ``almost connectedness'', i.e. we assume that $G$ has only a finite number of connected components. 

Many of the results of the paper \cite{kn:repre1}, can be summarized in the assertion displayed below that is valid under the hypothesis of the finiteness of the connected components mentioned above.  

\medskip

The following are equivalent: a) $\mathcal R(G)$ is finitely generated; b) $\pi^+(G^+)=A(G)$; c) if we call $K$ the clausure of the commutator subgroup of the connected component of the identity $G_1$, then $G/K$ is compact; d) if $\rho$ is a representation of $G$, then $\rho^+(G^+)$ is an affine algebraic group; e) $\pi(G_1)$ is the identity component of $B(G)$. 

Moreover, with respect to the natural structure of complex linear group of $G^+$ and assuming the validity of any of these additional hypothesis a)--d), the complex  representations of $G^+$ become algebraic.

\medskip

Concerning the second paper \cite{kn:repre2}, that could be broadly described as an investigation on the validity of the equivalent conditions a)--e) listed above in more general situations, we will only mention the following very neat result.

Assume that $h: G \rightarrow \C \in \operatorname{Hom}(G,\C)$ is a continuous multiplicative--additive morphism of $G$, and consider $\operatorname{exp}(h): G \rightarrow \C^\times$ its exponential, that clearly is a group like element of $\mathcal R(G)$.  Then, the authors prove that $\pi^+(G^+)=\{\alpha \in A(G): \alpha(\operatorname{exp}(h))= \operatorname{exp}(\alpha(h)), \forall \,\,h \in \operatorname{Hom}(G,\C)\}$ or equivalently if we write $\alpha'(f)= \alpha(f)(1)$ for $\alpha \in A(G)$ and $f \in \mathcal R(G)$, then $\pi^+(G^+)=\{\alpha \in A(G): \alpha'(\operatorname{exp}(h))= \operatorname{exp}(\alpha'(h)), \forall \,\, h \in \operatorname{Hom}(G,\C)\}$. It is clear that this condition is a generalization of the fact that condition c) above, implies condition b).

\medskip

In \cite{kn:repre3} that is the third paper of the series, the authors concentrate in the complex situation.

They define the concept of {\em reductive} complex analytic group:  a group is reductive if it has a faithful complex analytic representation and morever, all complex analytical representations are semisimple. 

In this paper the authors prove the analogous of the main theorem of the first paper.

Assuming that $G$ has a faithful complex analytic representation, then the following conditions are equivalent: a) The algebra of representative functions is finitely generated; b) 
$\pi(G)=A(G)$; c) The quotient of $G$ by its commutator subgroup is reductive; d) $G$ is the semidirect product of the radical of the commutator subgroup of $G$ and a reductive complex analytic group; e) $G$ can be identified with a complex linear group, where the complex analytic representations become the rational representations. 

\subsection{Representation and cohomology theory of groups with additional structure}

There was another direction of his early work that also led Hochschild to pay an increasing interest to the study of algebraic groups. 

Once that the homological techniques that permitted  to have a certain degree of control of the category of representations of finite groups were well established, it was  natural to apply the already developed methods to the study of representations of groups with additional structure. 

After his foundational work in homological algebra in general and in discrete group cohomology and its applications in particular, and armed with the techniques he had developed in his work in Lie and algebraic groups, Hochschild --usually with Mostow-- started to pay attention to these kind of situations. 

\medskip

We start by commenting three papers --Hochschild was the coauthor of two of them, where the authors discuss the problem  of {\em the extension of representations} in the situation of Lie groups and algebraic groups.  

The first two: \cite{kn:extlie1} and \cite{kn:extlie2} deal with the following problem in the category of Lie groups. 

\begin{defi}\label{defi:extension} Assume that $K \subset G$ is an inclusion in the category of Lie groups --or finite, analytic, affine algebraic, etc. A representation $V$ of $K$ extends to a representation $W$ of $G$, if $V \subset W$ as $K$ modules.  
\end{defi}

The basic questions the authors try to solve in both papers is if for an arbitrary pair of a group and a subgroup, every finite dimensional representation of $K$ admits a finite dimensional extension. 

In the case of finite groups, the problem has an obvious positive answer, given a representation $V$ of $K$ we cant take $W= \operatorname{Ind}_K^G(V)$, to be the corresponding induced $G$--module. 

In the case of Lie groups the situation is much more akward, considering that the natural construction given by induction even if it were available, it is not expected to produce a finite dimensional result.  The authors only obtain positive answers under the hypothesis that $K$ is normal in $G$ and with additional restrictions.  For example, the main result of \cite{kn:extlie1} reads as follows. 
\medskip

{\em Let $K \subset G$ be as above, with $K$ normal. Assume that there is an analytic subgroup $H$ of $G$ such that: $G=HK$, $H \cap K$ is compact, there is a finite dimensional representation of $H$ that is faithful on $H \cap K$.  Let $\rho$ be a representation of $K$, then $\rho$ can be extended to a representation of $G$ if and only if $\rho'([\operatorname{rad}(K), G])=1$, where the bracket represents the commutator subgroup, $\operatorname{rad}(K)$ is the radical of $K$ and $\rho'$ is the semisimple representation associated to $\rho$.}  Even though the formulation of the theorem is rather restrictive, it turns out to be very useful to simplify the proofs of some theorems concerning the existence of faithful representations due to E. Cartan, Got\^o, Malcev, etc. Part of the paper under consideration is devoted to present those proofs.

The situation treated in \cite{kn:extlie2}, is similar and refers to the case of an analytic group $G$ and a closed, normal subgroup $K$. The obtained results are also similar, the main difference being methodological as the author uses in a systematic way some of his own results in the theory of affine linear groups concerning linearly reductive subgroups and what is now called Mostow decomposition. This decomposition theorem guarantees that 
a connected affine algebraic group in characteristic zero, is the semidirect product of its unipotent radical and a linearly reductive group, see for example \cite{kn:libro} for a proof. 

Once the methods of the theory of affine algebraic groups were used to deal with the situation of analytic groups, it was clear that there were good perspectives that a reasonable theory of extensions could be developed in the algebraic environment. 

This line of work started with the publication of the joint paper with Białynicki--Birula and Mostow entitled: {\em Extensions of representations of algebraic linear groups}, see \cite{kn:obser}.

\medskip

Assume that $G$ is an affine algebraic group and $K$ a closed subgroup and that we are working over an algebraically closed field of arbitrary characteristic. 

\begin{defi} 
\begin{enumerate}
\item The subgroup $K$ is said to be observable in $G$ if the extension problem has a solution for all finite dimensional representations of $K$.
\item A character $\gamma: K \rightarrow \k$ is said to be extendible, if there exists a non zero regular function $f \in \k[G]$ with the property that for all $x \in K$ we have that $x\cdot f= \gamma(x)f$. 
\end{enumerate}
\end{defi} 

The main result proved in the paper can be summarized as follows.

\begin{teo} In the situation above the following are equivalent: a) The subgroup $K$ is observable in $G$; b) the homogeneous space is a quasi--affine variety; c) All the characters of $K$ are extendible to $G$.
\end{teo}

Even though Hochschild never published further in this direction, immediately after the publication of \cite{kn:obser}, 
the concept of observability started to be researched and rich applications were found; one specially remarkable is the application to invariant theory in the work of Grosshans and others. Recently, the concept of observable subgroup was generalized to the concept of observable action.  

The concept of extension of a representation considered above, has also been strengthened. 

In the notations above, a representation $V$ of $K$ is said to be strongly extendable if it is extendable and it has an extension $W$ with the additional property that the $G$--fixed part of $W$ coincides with the $K$--fixed part of $V$. 

In the 1970s, Cline, Parshall and Scott, defined the subgroup $K$ to be strongly observable in $G$ if all its representations are strongly extendable. It was proved that a subgroup $K$ is strongly observable if and only if it is exact --in the sense that the induction functor is exact-- and this is also equivalent to the property that the homogeneous space $G/K$ is an affine variety.  

Moreover, the above mentioned authors prove that in the case that $K$ is exact in $G$, then $\k[G]$ is a rationally injective $K$--module and conversely. This result settles for arbitrary characteristic the question of finding conditions that guaranteed that $\k[G]$ is cohomologically trivial as a $K$--module, that Hochschild asked in \cite{kn:ratinj} and that he solved in this paper for characteristic zero --Theorem 3.1 of \cite{kn:ratinj}. The question is important because its positive answer permits to have the machinery of the Hochschild--Serre spectral sequences, in working conditions. 

Next we consider Hochschild's work in the cohomology theory of affine algebraic groups. 

In the introduction to the paper {\em Cohomology of algebraic linear groups}, published in 1961 --see \cite{kn:calg}-- the author writes: ``{\em The theory of rational representations  of algebraic linear groups over fields of characteristic zero has, for some time, been in a sufficiently well developed state to call for an adaptation of homological algebra to the requisite category of `rational modules'\/''.}

In the intent of reproducing all the machinery of homological algebra to the context of the rational representations of an affine algebraic group, some important adaptations should be implemented.

One crucial point --that was noticed by Mostow in his paper {\em Cohomology of topological groups and solvmanifolds}, see \cite{kn:mostowcoho}, and developed in the paper we are considering-- is the following: contrary to the usual situation of discrete groups where the projective part of the machinery can be applied and it is somewhat more natural, in the case of Lie or algebraic groups one has to restrict the attention to injective resolutions. The point being that the categories under consideration have enough injectives but they need not have enough projectives. 

Once this machinery is available, one can define the rational cohomology groups of the affine algebraic group $G$ with coefficients in a rational representation $M$. We denote this family as: $\{H^i(G,M): i \geq 0\}$ and it can be defined explicitly in terms of the same cochain complex that is used to define group cohomology but with the additional restriction that all the functions are polynomial.

If we call $\mathfrak g$ the Lie algebra of $G$ and $E(G)$ the complex of differential forms based on the representative functions, i.e. on $\K[G]$, the main result proved in this paper is the following: \[H^\bullet(\mathfrak g, M)= H^\bullet(E(G)) \otimes H^\bullet(G,M).\] 

This isomorphism comes from the identification the  different elements of a particular instance of the Hochschild--Serre spectral sequence for Lie algebras as follows.

Call $G_u$ the unipotent radical of $G$, and  $\mathfrak g$ and $\mathfrak n$ the Lie algebras of $G$ and $G_u$ respectively and assume that $M$ is a rational $G$--module. 
  
From the spectral sequence for $\mathfrak g$, $\mathfrak n$ and $\mathfrak g/\mathfrak n$ and the fact that the last written Lie algebra is semisimple, one can prove the existence of an isomorphism $H^\bullet(\mathfrak g, M)= H^\bullet(\mathfrak g/\mathfrak n,\k) \otimes H^\bullet(\mathfrak n,M)^{\mathfrak g/\mathfrak n}$. By a direct identification of the right terms of the above isomorphism, the authors deduce the mentioned result.  It is important to mention that for the identification of  $H^\bullet(E(G))$ with $H^\bullet(\mathfrak g/\mathfrak n,\k)$
that comes from an extension of the dual of the canonical map $\mathfrak g \rightarrow \mathfrak g/\mathfrak n$ it is needed as a crucial ingredient the existence of Mostow semidirect product decomposition of $G=G_u \rtimes L$ for a linearly reductive $L$ that is only valid in characteristic zero.
One interesting consequence of the study of the cohomology or differential forms $H^\bullet(E(G))$ and the identification used above is the following: the differential form cohomology of the ring $\K[G]$ of polynomials in $G$ is trivial, if and only if $\K[G]$ is a polynomial ring. 

The line of work started in the paper just considered continues naturally with the joint paper with B. Kostant: {\em Differential forms and Lie algebra cohomology for algebraic linear groups}, that appeared one year later in 1962, see \cite{kn:dflacag}.  The purpose of this paper is to extend the results on differential forms to the situation of a homogeneous space of the form $G/H$ with $H$ a linearly reductive group --in characteristic zero. We use the same notation than before and call $\mathfrak h$ the Lie algebra of $H$. The isomorphism $H^\bullet(\mathfrak g, M)= H^\bullet(E(G)) \otimes H^\bullet(G,M)$ now becomes $H^\bullet(\mathfrak g, \mathfrak h, M)^H= H^\bullet(E(G/G_u \rtimes H)) \otimes H^\bullet(G,M)$, where the left terms indicates the relative Lie algebra cohomology. 

Moreover, the results just described were used in a later paper to produce an interesting cohomological characterization of the maximality of a linearly reductive subgroup of $G$, see \cite{kn:chahs}.

Along the way it is proved that for complex reductive
homogeneous spaces the de Rham cohomology can be computed using only holomorphic
differential forms, a result that was later vastly generalized by A. Grothendieck.

In a 1967 paper jointly written with Mostow and dedicated to Nakayama, this line of research is continued with the comparison of the rational cohomology with the holomorphic cohomology in case that $G$ is a complex analytic linear group, the basic result being --as expressed in the paper by the authors-- that: ``{\em the usual abundantly used connections between complex analytic representations of complex analytic groups and rational representations of algebraic groups extend[s] fully to the superestructure of cohomology''.} In fact they prove that the cohomology using holomorphic cochains is naturally isomorphic to the rational cohomology defined using polynomial cochains as in \cite{kn:calg}.

Tracing our steps back a little, we say some words about the important paper \cite{kn:coholie} published in 1962 jointly with Mostow and called: {\em Cohomology of Lie groups.}

This paper has the intention to compare for a Lie group $G$ and a finite dimensional $G$--module $V$, the different kind of cohomology groups with  coefficients in $V$ that could be defined using: continuous, differentiable or representable cochains --and eventually if the group is algebraic using polynomial cochains. 

The case of an affine algebraic group has been analyzed in \cite{kn:calg} and the case of representative cochains appears in this paper. The continuous cohomology appeared in \cite{kn:mostowcoho} but some parts of the theory developed there need to be reconsidered in order to adapt it to the machinery of injective resolutions --that as we mentioned before is the only way to work in the case of a group with superestructure. In all of the three cases mentioned above, one of the main achievements of the theory, is to link the global cohomology to the local cohomology of the Lie algebras. The two main theorems are the following: if $G$ is a real linear algebraic  group, $G_u$ its unipotent radical and $V$ a rational $G$--module, then there is an isomorphism $H^\bullet_{\text{rac}}(G,V) \otimes  H^\bullet_{\text{rep}}(G/G_u,\R) \cong H_{\text{rep}}^\bullet(G,V)$, and there is an analogous result for linear groups but one has to change, the representative cohomology by continuous cohomology and the rational cohomology by representative cohomology. 

\subsection{Unipotent groups in invariant theory}

We will now describe the article appearing in 1973 written together with Mostow, see \cite{kn:hmit}. Even though it is his only article dealing with ``hard core'' invariant theory, it represents a very important contribution as it was a breakthrough in a classical problem in the hardest case of the invariants of unipotent groups. It opened up what is today an active area of research with important open problems in train of being solved. 

The problem of the finite generation of rings of invariants, famously known as Hilbert's 14th problem, can be formulated as follows. 

{\em Let $V$ be a finite dimensional vector space and $H \subset \operatorname{GL}(V)$ a subgroup acting on $\k[V]$ with the induced linear action. Is the algebra of invariants $\k[V]^H$ finitely generated?}

As the condition of being an invariant is Zariski closed, one can take $H$ to be a closed subgroup of the linear group and then the problem can be attacked using the standard tools of the theory of affine algebraic groups. 

The important geometric meaning of the finite generation of the rings of invariants is clear.

Indeed, if we have an affine algebraic group $H$ acting on an affine variety with ring of polynomials $R$, as one expects that the invariant polynomials $R^H$ separate the geometric orbits --at least generically-- to have a finite number of algebra generators of $R^H$, will produce a finite criteria to decide whether two points of the variety are in the same orbit or not. In other words, the corresponding classification problem of the points of the variety can be solved in a finite number of steps by evaluation of a finite number of functions.

The search for groups with the above finiteness condition on invariants {\em for all actions} (in other words that are always adequate for taking quotients by them), culminated around mid 1970s thanks to the efforts of Haboush, Mumford, Nagata and Popov --to name those who in the opinion of this author were the main contributors, who showed that the only class of groups that guarantee the finiteness of the invariants, is the class of reductive groups.

In 1958, Nagata constructed a counterexample to Hilbert's problem, that consisted of a unipotent group of dimension 13 acting linearly on a space of dimension 32. Later, many counterexamples of smaller size have been found. This counterexample played a crucial role in the results of Haboush--Mumford--Nagata--Popov mentioned above. 

Moreover, once we know that the invariants of reductive groups on finitely generated algebras are finitely generated --another result of Nagata-- it is clear that the generic obstruction to finite generation is in the case of actions of unipotent groups. 

For that reason, the paper {\em Unipotent groups in invariant theory}, published in the Proceedings of the National Academy of Sciences, see \cite{kn:hmit},  dealing with the finite generation of the invariants of unipotent groups in $G$--module algebras, was considered extremely important. 

More precisely, the authors prove the following main result.
 
\begin{teo}\label{teo:hmit} Let $\K$ be an algebraically closed field of characteristic zero, $G$ a connected reductive group and $R$ a commutative finitely generated rational $G$--module algebra. If $U$ is a maximal unipotent subgroup of $G$, then the subalgebra $R\,^U = \{r \in R: u \cdot r = r, \quad \forall u \in U \} \subset R$ is finitely generated.  
\end{teo}

This theorem is one of the first general results dealing with invariants of unipotent groups, and can be interpreted as a generalization of the so called Weitzenbock's theorem\footnote{In accordance to A. Borel in \cite{kn:borelhistory}, Weitzenbock's theorem should better be called Maurer's theorem.} that guarantees the finite generation of the invariants of the additive group of a field acting linearly on a vector space. In that case $G$ is the special two by two linear group.

A particular case of Theorem \ref{teo:hmit} had been proved before by Dz. Hadziev in \cite{kn:rusosunipotentes}. It is proved for groups over the complex numbers and assumes that the commutative algebra $R$ is graded, and that the action of $G$ is homogeneous\footnote{The reader should take into account that to pass from the graded to the non graded case is not easy if we are not dealing with invariants of reductive groups.}.

Soon afterwords, F. Grosshans wrote a paper introducing an interesting viewpoing in invariant theory, that is related to the ideas of \cite{kn:hmit}. Building on the idea of observable subgroup and with a view of generalizing the situation treated by Maurer--Weitzenbock in dealing with invariants of the additive group of the field, the author defines $H \subset G$ to be what was later called a {\em Grosshans subgroup}, see \cite{kn:invgross}, and gives a very nice geometric characterization in terms of dimensions of the orbits of the action of the group $G$ on a certain representation. 

\begin{defi} Let $G$ be an affine algebraic group and $H \subset G$ a closed subgroup. The pair $(H,G)$ is called a Grosshans pair if: (1) $H$ is observable in $G$; (2) $\k[G]^H$ is a finitely generated algebra. 
\end{defi}

The observability condition is not too restrictive as we can always substitute the group $H$ by its observable closure in $G$ --see for example Chapter 12, Section 4, Lemma 5.2 of \cite{kn:libro}. 

The following isomorphism called ``the transfer principle''\footnote{It is also known as Grosshans principle, even if it was already known in particular cases to classical invariant theorists like Capelli or Roberts that used yet another name: ``{\em adjunction principle''.}} clarifies some problems related to invariants and justifies the above definition.

Assume that $H \subset G$ is a pair of affine algebraic groups and that $V$ is a rational $G$--module. Then there is a natural ismorphism between  $(\K[G]^H \otimes N)^G \cong N^H$.

Applying the transfer principle to the situation of Hilbert's problem and recalling that $\operatorname{GL}(V)$ is reductive, we deduce that that if $\k[\operatorname{GL}(V)]^H$ is finitely generated, then $\k[V]^H$ is also finitely generated. 

More generally, if $(H,G)$ is a Grosshans pair and $G$ is reductive, we deduce that the $H$--invariants of a rational affine $G$--module algebra are finitely generated. 

The main theorem of \cite{kn:hmit} can be formulated as follows: if the base field has characteristic zero, $G$ is a reductive group and $U$ a maximal unipotent subgroup, then $(U,G)$ is a Grosshans pair.

A far ranging conjecture that as far as the author knows is still open, and that in case it is true is an ample generalization of the seminal result by Hochschild and Mostow is the following.  
 
Conjecture of Popov--Pommerening: {\em If $G$ is a reductive group and $U \subset G$ is a unipotent subgroup normalized by a maximal torus, then $(U,G)$ is a Grosshans pair.}

\medskip


\end{document}